\documentclass[12pt]{amsart}
\usepackage{amsmath,amssymb}
\newtheorem{theorem}{Theorem}
\newtheorem{proposition}[theorem]{Proposition}
\newtheorem{corollary}[theorem]{Corollary}

\newtheorem{example}[theorem]{Example}
\def\C{\Bbb C}
\def\D{\Bbb D}
\def\R{\Bbb R}
\def\O{\mathcal O}
\def\eps{\varepsilon}
\def\om{\omega}
\def\g{\gamma}
\def\c{\tilde c}

\def\z{\tilde z}
\def\hz{\hat z}
\def\vz{\check z}
\def\w{\tilde w}
\def\hw{\hat w}
\def\vw{\check w}

\def\ds{\displaystyle}

\title{Estimates of the Kobayashi and quasi-hyperbolic distances}

\author{Nikolai Nikolov and Lyubomir Andreev}

\address{N. Nikolov: Institute of Mathematics and Informatics\\Bulgarian Academy
of Sciences\\Acad. G. Bonchev 8, 1113 Sofia, Bulgaria\newline
\indent Faculty of Information Sciences\\
State University of Library Studies and Information Technologies\\
Shipchenski prohod 69A, 1574 Sofia, Bulgaria}\email{nik@math.bas.bg}

\address{L. Andreev: Institute of Mathematics and Informatics\\Bulgarian Academy
of Sciences\\Acad. G. Bonchev 8, 1113 Sofia, Bulgaria}
\email{lyubomir.andreev@math.bas.bg}

\subjclass[2010]{32F45, 51M10}

\keywords{Kobayashi distance, quasi-hyperbolic distance}

\thanks{L. Andreev is partially supported by the Bulgarian national science found
under contract DFNI-I 02/14.}

\begin{document}

\begin{abstract}{Universal upper bounds for the Kobayashi and quasi-hyper\-bolic distances
near Dini-smooth boundary points of domains in $\C^n$ and $\R^n,$ respectively,
are obtained.}
\end{abstract}

\maketitle

\section{Introduction and results}

Let $D$ be a domain in $\C^n.$ The Kobayashi (pseudo)distance
$k_D$ is the largest (pseudo)distance not exceeding the Lempert function
$$l_D(z,w)=\inf\{\tanh^{-1}|\alpha|:\exists\varphi\in\O(\D,D)
\hbox{ with }\varphi(0)=z,\varphi(\alpha)=w\},$$ where $\D$ is the
unit disc. Note that $k_D$ is  the integrated form of the
Kobayashi (pseudo)metric
$$\kappa_D(z;X)=\inf\{|\alpha|:\exists\varphi\in\O(\D,D) \hbox{
with }\varphi(0)=z,\alpha\varphi'(0)=X\}.$$
We also point out that $k_D=l_D$ for any domain $D$ in $\C.$

Various estimates of invariant distances and metrics play crucial role in many questions
in complex analysis.

Set $d_D(z)=\mbox{dist}(z,\partial D).$ Recall the following upper bound for $k_D.$

\begin{proposition}\label{FR}\cite[Proposition 2.5]{FR} Let $a$ be a $\mathcal C^{1+\eps}$-smooth
boundary point of a domain $D$ in $\C^n.$ Then there exist a neighborhood $U$ of $a$
and a constant $c>0$ such that for $z,w\in D\cap U,$
$$2k_D(z,w)<\log\left(1+\frac{|z-w|}{d_D(z)}\right)+\log\left(1+\frac{|z-w|}{d_D(w)}\right)+c.$$
\end{proposition}

This easily implies the next corollary.

\begin{corollary} Let $D$ be a $\mathcal C^{1+\eps}$-smooth bounded domain
in $\C^n.$ Then there exists a constant $c>0$ such that
$$2k_D(z,w)<\log\left(1+\frac{|z-w|}{d_D(z)}\right)+\log\left(1+\frac{|z-w|}{d_D(w)}\right)+c.$$
\end{corollary}

The following weaker inequality holds for $l_D.$

\begin{proposition}\cite[Theorem 1]{NPT} Let $D$ be a $\mathcal C^{1+\eps}$-smooth
bounded domain $D$ in $\C^n.$ Then there exists a constant $c>0$ such that
$$2l_D(z,w)<-\log d_D(z)-\log d_D(w)+c.$$
\end{proposition}

Note that $C^1$-smoothness is not enough for the last two inequalities.

\begin{example}\label{exa}\cite[Example 2]{NPT} The image of $\D$
under the (univalent) map $z\mapsto 2z+(1-z)\log(1-z)$ is a $\mathcal C^1$-smooth
bounded domain $D$ for which
$$\lim_{\R\ni z\uparrow 2}(2k_D(z,0)+\log d_D(z))=+\infty.$$
\end{example}

Recall now that the simple inequality $\ds\kappa_D(z;X)\le\frac{|X|}{d_D(z)}$ shows that
the Kobayashi distance does not exceed the quasi-hyperbolic distance
$$h_D(z,w)=\inf_{\g}\int_\g\frac{|du|}{d_D(u)},$$
where the infimum is taken over all rectifiable curves $\g$ in $D$ joining $z$ to $w$
(by \cite[Lemma 1]{GO}, the infimimum is attained). The definition of $h_D$ in the real case
is the same. This distance arises in the the theory of quasi-conformal maps.

We have more if $\mathcal C^1$-smoothness is required (despite Example \ref{exa}).

\begin{proposition}\label{sup} (a) If $a$ is a $\mathcal C^1$-smooth
boundary point of a domain $D$ in $\C^n,$ then
$$\limsup_{\substack{z,w\to a\\z\neq w}}\frac{k_D(z,w)}{h_D(z,w)}\le\frac12.$$

\noindent (b) If $D$ is a $\mathcal C^1$-smooth bounded domain in $\C^n,$ then
$$\limsup_{z\to\partial D}\frac{k_D(z,w)}{h_D(z,w)}\le\frac12\quad\mbox{ uniformly in }w\neq z.$$
\end{proposition}

Set $\ds s_D(z,w)=2\sinh^{-1}\frac{|z-w|}{2\sqrt{d_D(z)d_D(w)}}.$ Note that $h_D=s_D$ if $D$
is a half-space in $\R^n$ (cf. \cite[(2.8)]{Vuo}).

\begin{proposition}\label{lim} (a) If $a$ is a $\mathcal C^1$-smooth
boundary point of a domain $D$ in $\R^n,$ then
$$\lim_{\substack{z,w\to a\\z\neq w}}\frac{h_D(z,w)}{s_D(z,w)}=1.$$

\noindent (b) If $D$ is a $\mathcal C^1$-smooth bounded domain in $\R^n,$ then
$$\lim_{z\to\partial D}\frac{h_D(z,w)}{s_D(z,w)}=1\quad\mbox{ uniformly in }w\neq z.$$
\end{proposition}

Let $\ds q_D(z,w)=\frac{h_D(z,w)}{s_D(z,w)}$ if $z,w\in D,$ $z\neq w,$ and $q_D(z,w)=1$
otherwise. Since $\ds\lim_{\substack{z,w\to u\\z\neq w}}q_D(z,w)=1$ for any $u\in D,$
Proposition \ref{lim} (b) means that $q_D$ is a continuous function on $\R^{2n}.$

Our main purpose is, refining the proof of Proposition \ref{FR}, to extend
this proposition in two directions. First, we put the constant $c>0$ in
front of the factor $|z-w|$ which agrees to the case $z=w.$ Moreover, we choose
$c$ to be universal. Second, we replace the $\mathcal C^{1+\eps}$-smoothness
by the weaker assumption of Dini-smoothness.

Recall that a $\mathcal C^1$-smooth boundary point $a$ of a domain $D$ in $\C^n$ is said to be
Dini smooth (or Lyapunov-Dini smooth) if the inner unit normal vector $n$ to $\partial D$ near $a$
is a Dini-continuous function. This means that there exists a neighborhood $U$ of $a$ such that
$\ds\int_0^1\frac{\om(t)}{t}dt<+\infty,$ where
$$\om(t)=\om(n,\partial D\cap U,t):=\sup\{|n_x-n_y|:|x-y|<t,\ x,y\in \partial D\cap U\}$$
is the respective modulus of continuity.

For $c>0$ set
$$v^c_D(z,w)=2\log\left(1+\frac{c|z-w|}{\sqrt{d_D(z)d_D(w)}}\right).$$
By \cite[Theorem 1.1]{DHV}, $v^c_D$ is a distance if $c\ge 2.$ Let $\ds c_0=1+\frac{\sqrt 2}{2}.$

\begin{theorem}\label{main} Let $a$ be a Dini-smooth boundary point of a domain $D$ in $\C^n,$
resp. $\R^n.$ Then for any constant $c>c_0$ there exists a neighborhood $U$ of $a$ such that
on $(D\cap U)^2,$
$$2k_D\le v^c_D,\mbox{ resp. }h_D\le v^c_D.$$
\end{theorem}

\begin{corollary}\label{cor} Let $D$ be a Dini-smooth bounded domain
in $\C^n,$ resp. $\R^n.$ Then there exists a constant $c>c_0$ such that
$$2k_D\le v^c_D,\mbox{ resp. }h_D\le v^c_D.$$
\end{corollary}

Note now that $h_D$ has a minorant $i_D$ for which $v^{1/2}_D\le i_D\le s_D\le v^1_D.$

\begin{proposition}\label{GHM}\cite[Lemma 2.6]{GHM} If $D$ is a proper subdomain of $\R^n,$  then
$$h_D(z,w)\ge i_D(z,w):=2\log\frac{d_D(z)+d_D(w)+|z-w|}{2\sqrt{d_D(z)d_D(w)}}.$$
\end{proposition}

Observe that $h_D=i_D$ if $D\subset\R$ (then $D$ is an open interval or ray).

The next corollary is a consequence of Theorem \ref{main}, Corollary \ref{cor} and
Proposition \ref{GHM}.

\begin{corollary} If $a$ is a Dini-smooth boundary point of a domain $D$ in $\C^n,$
then
$$\limsup_{z,w\to a}(2k_D(z,w)-h_D(z,w))\le2\log(2+\sqrt2).$$

Therefore, if $D$ is a Dini-smooth bounded domain in $\C^n,$
then the function $2k_D-h_D$ is bounded from above.
\end{corollary}

The next proposition shows that $i_D$ is a distance. Then the equality
$\ds\lim_{\substack{z,w\to u\\z\neq w}}\frac{i_D(z,w)}{|z-w|}=\frac{1}{d_D(u)}$
implies that $h_D$ is the intrinsic distance of $i_D;$ in particular, $h_D\ge i_D.$

\begin{proposition}\label{dist} If $X$ is a metric space and $f:X\to\R^+$ is
an 1-Lip\-schitz function, then
$\ds\rho(z,w)=\log\frac{f(z)+f(w)+|z-w|}{2\sqrt{f(z)f(w)}}$ is a distance.
\end{proposition}

The rest of the paper is organized as follows. Section 2 contains the proof
of Propositions \ref{sup}, \ref{lim} and \ref{dist}. Section 3
contains the proof of Theorem \ref{main}. Section 4 is devoted to an application
of Corollary \ref{cor} to boundary continuity of proper holomorphic maps.

\section{Proofs of Propositions \ref{sup}, \ref{lim} and \ref{dist}}

\noindent{\it Proof of Proposition \ref{sup}.} We shall use that
if $a\in\partial D$ is  $\mathcal C^1$-smooth, then (cf. \cite[(2)]{N1})
\begin{equation}\label{ine}
\limsup_{z\to a}\max_{|X|=1}\kappa_D(z;X)d_D(z)\le\frac12.
\end{equation}

First, we shall prove (b). By \eqref{ine}, for any $c>1/2$
there exists a neighborhood $U$ of $\partial D$ such that
$$2\kappa_D(z;X)d_D(z)\le c|X|,\quad z\in D\cap U.$$

Consider a geodesic $\g$ of $h_D$ joining $z\in D\cap U$ and $w\in D.$ If $\g\in U,$ then
$k_D(z,w)\le c h_D(z,w).$ Otherwise, let $u$ and $v$ be the first and the last points,
where $\g$ meets $\partial U$ (possible $u=v$). Then
$$k_D(z,w)\le ch_D(z,w)+k_D(u,v).$$
Since
$\ds\lim_{z\to\partial D}\inf_{u\in D\setminus U}h_D(z,u)=+\infty$
by Proposition \ref{GHM}, it follows that
$$\limsup_{z\to\partial D}\frac{k_D(z,w)}{h_D(z,w)}
\le c\quad\mbox{ uniformly in }w\neq z.$$
It remains to let $c\to 1/2.$

The proof of (a) is the same, choosing a respective neighborhood $U$ of $a.$
\smallskip

\noindent{\it Proof of Proposition \ref{lim}.} (a) After translation and rotation,
we may assume that $a=0$ and that there is a neighborhood $U$ of $0$ such that
$$D':=D\cap U=\{x\in U: r(x):=x_1+f(x')>0\},$$
where $f$ is a $\mathcal C^1$-smooth function in $\R^n$ with $f(0)=0$ and $\nabla f(0)=0.$

Let $c>1$ and $\theta(x)=(r(x),x').$ We may shrink $U$ such that
\begin{equation}\label{c}
c^{-1}|x-y|\le |\theta(x)-\theta(y)|\le c|x-y|,\quad x,y\in U.
\end{equation}

Choose now a neighborhood $V\subset U$ of $0$ such that $d_{D'}=d_D$ on $D\cap V.$
Using, for example, \cite[Theorem 2]{GO}, one can find a neighborhood
$W\subset V$ of $0$ such that any geodesic for $h_D(z,w)$ belongs to $D\cap V$ if
$z,w\in \tilde D=D\cap W.$ Then $h_D=h_{D'}$ on $\tilde D^2.$

Set $\Pi=\{x\in\R^n:x_1>0\}$ and $\Pi'=\theta(D').$ Using the above arguments,
we may shrink $W$ such that $h_\Pi=h_{\Pi'}$ on $(\theta(\tilde D))^2.$

On the other hand, \eqref{c} implies that (cf. \cite[Exercise 3.17]{Vuo})
$$c^{-2}h_{D'}(z,w)\le h_{\Pi'}(\theta(z),\theta(w))\le c^2h_{D'}(z,w),\quad z,w\in D'.$$

Let $z,w\in \tilde D.$ Then
$$c^{-2}h_D(z,w)\le h_\Pi(\theta(z),\theta(w))\le c^2h_D(z,w).$$
Using \eqref{c} again, we get that
\begin{align*}h_\Pi(\theta(z),\theta(w))&=2\sinh^{-1}\frac{|\theta(z)-\theta(w)|}{2\sqrt{r_D(z)r_D(w)}}\\
&\le2\sinh^{-1}\frac{c^2|z-w|}{2\sqrt{d_D(z)d_D(w)}}\le c^2s_D(z,w).
\end{align*}
We obtain in the same way that
$$h_\Pi(\theta(z),\theta(w))\ge c^{-2}s_D(z,w).$$
So
$$c^{-4}h_D(z,w)\le s_D(z,w)\le c^4h_D(z,w)$$
which implies the desired result.
\smallskip

(b) It is an easy consequence of (a) and we omit the details.
\smallskip

\noindent{\it Proof of Proposition \ref{dist}.}
$\rho(z,u)+\rho(u,w)\ge \rho(z,w)\Leftrightarrow$
$$\frac{f(z)+f(u)+|z-u|}{2\sqrt{f(z)f(u)}}\cdot
\frac{f(u)+f(w)+|u-w|}{2\sqrt{f(u)f(w)}}
\ge\frac{f(z)+f(w)+|z-w|}{2\sqrt{f(z)f(w)}}$$
$$\Leftrightarrow(f(z)+f(u)+|z-u|)(f(u)+f(w)+|u-w|)$$
$$\ge 2f(u)(f(z)+f(w)+|z-w|)$$
$$\Leftrightarrow(f(z)-f(u)+|z-u|)(f(w)-f(u)+|w-u|)$$
$$+2f(u)(|z-u|+|u-w|-|z-w|)\ge 0.$$
The last inequality holds, since all factors are nonnegative.

\section{Proof of Theorem \ref{main}} First, we shall prove
the result for $k_D.$

We may find a neighborhood $V$ of $a$ and a number $s>0$ such that:

\noindent-- $\om=\om(n,\partial D\cap V,\cdot)$ is a Dini-continuous function;

\noindent-- if $b\in\partial D\cap V$ and $\phi_b(\tau)=b+\tau n_b$ ($\tau\in\C$),
then $\phi_b(3s\D)\subset V,$ and

$D_b:=\phi^{-1}_b(D\cap V)=\{\tau=\zeta+i\eta\in\phi^{-1}_b(V):\zeta+f_b(\eta)<0\};$

\noindent-- $\om_b\le\om,$ where $\om_b$ is the modulus of continuity
of $f'_b$ on $(-2s,2s).$

Note that $\om_\ast=\sup\om_b$ is an subadditive function. Then its modulus
of continuity does not exceed $\om$ and hence $\om_\ast$ is
Dini-continuous. So, we may find a Dini-smooth domain $G,$ symmetric with respect
to the real axis, such that
$$D_\ast\cap s\D\subset G\subset D_\ast\cap 2s\D\subset D_b,$$
where $\ds D_\ast=\{\tau\in\Bbb C:\zeta+\int_0^{|\eta|}w_\ast(t)dt<0\}.$

Note that there a conformal map $\theta:G\to\D$ which extends to a
$C^1$-diffeomorphism from $\overline G$ to $\overline \D$ (cf. \cite[Theorem
3.5]{Pom}). We may assume that $\theta(G\cap\R)=(-1,1).$

Let now $2c_1\ge q>1.$ Choose a ball $U$
of radius $\ds r\le\frac{s}{2c_1q^2+1}$ centered at $a$ such that if
$z\in D\cap U,$ then $d_D(z)=|z-z'|$ for some
$z'\in\partial D\cap V.$ For $z,w\in D\cap U$ set
$$\z=z+c_1q^2|z-w|n_{z'},\ \hz=d_D(z),\ \vz=d_D(z)+c_1q^2|z-w|.$$
Then $\hz,\vz\in G$ and $\phi_{z'}(\hz)=z,$ $\phi_{z'}(\vz)=\z.$
Hence
$$k_D(z,w)-k_D(\z,\w)\le k_D(z,\z)+k_D(w,\w)$$
$$\le k_G(\hz,\vz)+k_G(\hw,\vw)\le
k_\D(\theta(\hz),\theta(\vz))+k_\D(\theta(\hw),\theta(\vw)).$$

Since
$$k_\D(\alpha,\beta)=\frac12\log\left(1+\frac{2(\alpha-\beta)}
{(1-\alpha)(1+\beta)}\right),\quad -1<\beta\le\alpha<1,$$
$$\lim_{\nu\to 0}\frac{d_\D(\theta(\nu))}{d_G(\nu)}=|\theta'(0)|,\quad
\lim_{u\to a}\frac{d_G(\hat u)}{d_D(u)}=1,$$
we may shrink $r$ such that
\begin{equation}\label{cr}
k_D(z,w)\le k_D(\z,\w)+\frac12\log\left(\hskip-1mm1+\frac{c_1q^3|z-w|}{d_D(z)}\right)
\hskip-1mm\left(\hskip-1mm1+\frac{c_1q^3|z-w|}{d_D(z)}\right).
\end{equation}

To estimate $k_D(\z,\w),$ we shrink $r$ once more such that
$$|n_z-n_w|\le q-1,\ |\z-z'|\le q d_D(\z),\ |\w-w'|\le q d_D(\w).$$
Then
$$|\z-\w|\le q|z-w|<c_1^{-1}\min\{d_D(\z),d_D(\w)\}\le d_D(\z)+d_D(\w).$$
Choosing $V$ convex, it follows that $[\z,\w]\subset V.$ By \eqref{ine},
we may shrink $V$ such that
$$2\kappa_D(u;X)d_D(u)\le q|X|,\quad u\in D\cap V.$$

Now, we shall proceed as in the proof of \cite[Lemma 2.1]{GHM}.
There is a point $\vartheta$ in the complex line through $\z$ and $\w$
such that $|\vartheta-\z|=d_D(\z)$ and $|\vartheta-\w|=d_D(\w).$ If
$\gamma=[z,w],$ then
$$\frac{\tanh k_D(\z,\w)}{q}\le\tanh\frac{k_D(\z,\w)}{q}\le\tanh\int_\gamma\frac{|du|}{2d_D(u)}$$
$$\le\tanh\int_\gamma\frac{|du|}{2|u-\vartheta|}
=\frac{|\z-\w|}{d_D(\z)+d_D(\w)}\le\frac{q^2|z-w|}{|\z-z'|+|\w-w'|}$$
$$=\frac{q^2|z-w|}{d_D(z)+d_D(w)+2c_1q^2|z-w|}\le
\frac{q^3|z-w|}{d_D(z)+d_D(w)+2c_1q^3|z-w|}.$$

It follows from here and \eqref{cr} that
$$2k_D(z,w)\le p_D(z,w,c_2):=\log\frac{g(c_2,\gamma,\delta)}{\gamma\delta},$$
where $c_2=c_1/q,$ $\ds\gamma=\frac{d_D(z)}{q^4|z-w|},$
$\ds\delta=\frac{d_D(w)}{q^4|z-w|}$ ($z\neq w$) and
$$g(\gamma,\delta,c_2)=(c_2+\gamma)(c_2+\delta)\frac{2c_2+1+\gamma+\delta}
{2c_2-1+\gamma+\delta}.$$

Assume for a while that $p_D(z,w,c_2)\le v^c_D(z,w),$ where $c=\c_0q^3.$
This means that $g(\gamma,\delta,c_2)\le(\sqrt{\gamma\delta}+\c_0)^2.$
Letting $\gamma,\delta\to 0,$ we obtain
$$c_2\sqrt\frac{2c_2+1}{2c_2-1}\le\c_0.$$
Letting $\gamma\to 0,q\to 1$ with $\delta=\gamma+q^4,$ we get
$c_2+1\le\c_0.$ Thus  
$$\c_0\ge\min_{c_2>1/2}\max\{c_2\sqrt\frac{2c_2+1}{2c_2-1},c_2+1\}=c_0$$
and $\c_0=c_0\Leftrightarrow c_2=\sqrt2/2.$

Note now that 
$$g(\gamma,\gamma+1,c_2)-g(\gamma,\delta,c_2)=\frac{\gamma+1-\delta}{2c_2-1+\gamma+\delta}
(2c_2^2-1+3c_2\gamma+c_2\delta+\gamma^2+\gamma\delta).$$

Let $c_2\ge\sqrt2/2.$ Since $\delta<\gamma+1,$ it follows that
$$g(\gamma,\delta,c_2)<g(\gamma,\gamma+1,c_2)=(\gamma+c_2+1)^2.$$

So, $2k_D\le v^c_D$ on $(D\cap U)^2,$ where $c=c_0q^3.$
\smallskip

The result for $h_D$ can be obtained in the same way.
If $n=1,$ then $h_D=i_D\le v^1_D.$

Let $n\ge 2.$ We may choose $G$ to be convex which implies that
$2d_G\kappa_G(\cdot;1)\ge1$ and hence $2k_G\ge h_G.$
Then we may proceed as above.

\section{An application of Corollary \ref{cor}}

The next result is inspired by \cite[Theorem 1.1]{FR} and \cite[Theorem 4.1]{Zim}.

\begin{proposition} Let $D$ be a Dini-smooth bounded domains in $\C^n$
with a negative plurisubharmonic function $u_D\ge-d_D$.\footnote{If $D$ is convex,
then $u_D=-d_D$ works.} Let $G$ be a convex Dini-smooth bounded
domains in $\C^n$ such that $\partial G\cap T_b^\C\partial G=\{b\}$ for any
$b\in\partial G.$ Then every proper holomorphic map $f:D\to G$ extends to a continuous
map from $\overline D$ onto $\overline G.$
\end{proposition}

\noindent{\it Proof.} Assume the contrary. Then
one may find sequences $(z^k_j)_j\subset D$ ($k=1,2$) such that for $w^k_j=f(z^k_j),$
$$\lim_{j\to\infty}z^k_j=a\in\partial D,\ \lim_{j\to\infty}w^k_j=b^k\in\partial G,
\ b^1\neq b^2,\ w^1_j\neq w^2_j.$$ By \cite[Lemma 3.5]{Zim} and Corollary \ref{cor},
there is a constant $c>c_0$ for which
$$-\log(cd_G(w^1_j)d_G(w^2_j))<2k_G(w^1_j,w^2_j)
\le 2k_D(z^1_j,z^2_j)\le v^c_D(z^1_j,z^2_j).$$

On the other hand, $v(w)=\max\{u(z):z\in f^{-1}(w)\}$ is a plurisubharmonic function on
$G.$ Since the Hopf lemma holds on $G$ (cf. \cite[Theorem 3.5]{Wid}), we may increase $c$
such that $d_G\circ f<-\sqrt cv\circ f\le \sqrt cd_D.$ Then
$$1<c\sqrt{d_D(z^1_j)d_D(z^2_j)}+c^2|z^1_j-z^2_j|$$
which is impossible for any $j$ large enough.

\end{document}